\documentclass[12pt]{amsart}
\def\R{\mathbb R}
\def\P{\mathbb P}
\def\C{\mathbb C}
\makeatletter
\let\over\@@over
\makeatother

\newtheorem{theorem}{Theorem}[section]
\newtheorem{lemma}{Lemma}[section]

\newtheorem{definition}{Definition}[section]

\newtheorem{remark}{Remark}[section]
\begin{document}
\title[Gauss curvature of certain minimal surfaces in $\R^m$]
{An estimate for the Gauss curvature of minimal surfaces 
in ${\R }^m$ whose Gauss map omits a set of hyperplanes}

\author{Robert Osserman}
\address{\hskip-\parindent
Robert Osserman
Department of Mathematics, Stanford University, Stanford, CA 94305;
 Mathematical Sciences Research Institute, 1000 Centennial Drive, 
Berkeley, CA 94720.}
\email{ro@msri.org}

\author{Min Ru}
\address{\hskip-\parindent
Min Ru,
Department of Mathematics, University of Houston, Houston, TX 77204-3476; 
 Mathematical Sciences Research Institute, 1000 Centennial Drive, 
Berkeley, CA 94720. }
\email{minru@math.uh.edu, minru@msri.org}

\thanks{Research supported in part by NSF grant DMS-9506424,
and at MSRI by NSF grant DMS-9022140.}

\begin{abstract}
We give an estimate of the Gauss curvature for 
minimal surfaces in ${\R}^m$ whose Gauss map omits more than 
$m(m+1)/2$ hyperplanes in ${\P}^{m-1}({\C})$.
\end{abstract}
\maketitle

\section{Introduction}
The purpose of this paper is to prove the following theorem.

\begin{theorem}[Main Theorem]  
Let $x: M \rightarrow {\R }^{m}$ be a minimal surface immersed in 
${\R }^m$. Suppose that its generalized Gauss map $g$ omits 
more than ${m(m+1)\over 2}$ hyperplanes in ${\P }^{m-1}({\C })$, 
located in general position. Then there exists a constant $C$, depending on 
the set of omitted hyperplanes, but not the surface, such that 
\begin{eqnarray} 
 |K(p)|^{1/2} d(p) \leq C 
\end{eqnarray}
where $K(p)$ is the Gauss curvature of the surface at $p$, 
 and $d(p)$ is the geodesic distance 
from $p$ to the boundary of $M$.
\end{theorem}
This theorem provides a considerable sharpening of an earlier result of 
the same type:
\begin{theorem}(Osserman [O1]) An inequality of the form (1) holds for all 
minimal surfaces in ${\R }^m$ whose Gauss
 map omits a neighborhood of some hyperplane in ${\P }^{m-1}({\C })$.
\end{theorem}

Also, Theorem 1.1 implies the earlier result:
\begin{theorem}(Ru[1])
Let $x: M \rightarrow {\R }^m$ be a complete minimal surface immersed in 
${\R }^m$. Suppose that its generalized Gauss map $g$ omits 
more than ${m(m+1)\over 2}$ hyperplanes in ${\P }^{m-1}({\C })$, 
located in general position. Then $g$ is constant and 
the minimal surface must be a plane. 
\end{theorem}

In fact, given any point $p$ on a complete surface satisfying the hypotheses,
inequality (1) must hold with $d(p)$ arbitrarily large, so that $K(p)= 0$. But
a minimal surface in ${\R }^m$ with $K \equiv 0$ must lie on a plane (see 
[HO1]) and hence its Gauss map $g$ is constant. 

\bigskip
Theorem 1.3 had been proved earlier by Fujimoto [F2]  in the case where 
the Gauss map $g$ was assumed nondegenerate. Fujimoto (see [F3])
also showed that the 
number $m(m+1)/2$ was optimal in that for every odd dimension $m$, there 
exist complete minimal surfaces whose Gauss map omits $m(m+1)/2$ hyperplanes 
in general position. It follows that Theorem 1.1 is also an 
optimal result of its type, since with any smaller number of omitted hyperplanes, a universal 
inequality  of the form (1) cannot be valid, at least in odd dimensions. 

\bigskip
When $m=3$, we may consider the classical Gauss map into the unit sphere. 
Fujimoto [F1] showed that an inequality of type (1) holds whenever the Gauss 
map omits 5 given points. In fact he obtained an explicit expression for $C$ 
in that case, depending on the given points. Ros [Ro1] gave a different proof 
that does not yield an explicit value for the constant $C$, but which allows the extension to higher dimension that we give here.

\section{Some theorems and lemmas}
In this section, we recall some results which will be used later.

We first recall the following construction theorem of minimal surfaces.
\begin{theorem}(see [CO1])
 Let $M$ be an open Riemann surface and let $\omega_1, \omega_2, \dots, \omega_m$ be holomorphic forms on $M$ having
 no common zero, no real periods and locally satisfying the 
identity
\[ f_1^2 + f_2^2 + \dots + f_m^2 = 0\]
for holomorphic functions $f_i$ with $\omega_i = f_i dz$. Set
\[ x_i = 2Re \int_{z_0}^z \omega_i,\]
for an arbitrary fixed point $z_0$ of $M$. Then the surface 
$x = (x_1, \dots, x_m): M \rightarrow {\R }^m$ is a minimal surface 
immersed in ${\R }^m$ such that the Gauss map is the map 
$g = [\omega_1: \dots :\omega_m]: M \rightarrow Q_{m-2}({\C })$ and 
the induced metric is given by 
\[ ds^2 = 2(|\omega_1|^2 + \dots +|\omega_m|^2).\]
\end{theorem}

\bigskip
The  following is the  general version of Hurwitz's theorem:

\begin{theorem}[Hurwitz's theorem] Let $f_j: M \rightarrow N$ be a sequence of 
holomorphic maps between two connected complex manifolds  
converging uniformly on every compact subset of $M$ to 
a holomorphic map $f$. If the image of each map $f_j$ misses a divisor 
$D$ of $N$ then either the image of $f$ misses $D$ or it lies entirely in $D$. 
\end{theorem}
Proof. Assume first that $D = \{z| g(z) = 0\}$ for some holomorphic function $g$. Then $g\circ f_j$ is a sequence of holomorphic functions converging 
to the holomorphic function $g \circ f$. Since $g \circ f_j$ is non-vanishing, 
by the classical Hurwitz theorem the limit function is either identically zero or non-vanishing. In other 
words the image of $f$ either lies entirely in $D$ or misses $D$ completely. 

In the general case, if $f$ does not miss $D$ entirely then there exists a 
point $q$ in $D$ and a point $p$ in $M$ such that $f(p) = q$. There exists a 
neighborhood $U$ of $q$ and a holomorphic function $g$ on $U$ so that 
$D \cap U = \{ z| g(z) = 0\}.$  Applying the previous argument to the restriction of the sequence of maps to the open set $V = f^{-1}(U)$ in $U$, we conclude 
that $f(V)$ is contained in $D \cap U$. Since $M$ is connected, the principle of analytic continuation implies that the image $f(M)$ is contained in $D$.

\bigskip
\begin{lemma}
Let $D_r$ be the disk of radius $r, 0 < r < 1$,
and let $R$ be the hyperbolic radius of $D_r$. Let 
\[ ds^2 = \lambda (z)^2 |dz|^2\]
be any conformal metric on $D_r$ with the property that the geodesic distance 
from $z=0$ to $|z| = r$ is greater than or equal to $R$. If the Gauss 
curvature $K$ of the metric $ds^2$ satisfies 
\[ -1 \leq K \leq 0\]
then the distance of any point to the origin in the metric $ds^2$ is greater 
than or equal to the hyperbolic distance. 
\end{lemma}

\begin{remark}
The hyperbolic metric in the unit disk is given by 
\[ d{\hat s}^2 = {\hat \lambda}(z)^2 |dz|^2, ~~{\hat \lambda}(z) = {2\over 1-|z|^2}, \]
and has curvature ${\hat K} \equiv -1$. The relation between the quantities $R$ and $r$ is therefore given by 
\[ R = \int_0^r {\hat \lambda}(z) |dz| = \int_0^r {2\over 1-t^2}dt = \log {1+r\over 1-r}\]
and the conclusion of the theorem is that 
\[ \rho (z) \ge {\hat \rho}(z) = \log {1+|z|\over 1 - |z|} \]
where $\rho$ and ${\hat \rho}$ represent the distances from 
the point $z$ to the origin in the metric $ds^2$ and the hyperbolic metric, 
respectively. 
\end{remark}
\begin{remark} Lemma 2.1 and its proof are basically geometric 
reformulations of lemma 6 of Ros[Ro 1]. The lemma may be viewed as a kind of dual to the Ahlfors 
form of the Schwarz-Pick  lemma [A 1].   
\end{remark}

\bigskip
\noindent{{\sl Proof of lemma 2.1}.}
Note first, that in the relation above between $R$ and $r$, we have 
\[ {dR\over dr} = {2\over 1-r^2} > 0\]
and we may solve for $r$ in terms of $R$:
\begin{eqnarray}
 r = {e^R-1\over e^R + 1},
\end{eqnarray}
or in general 
\begin{eqnarray}
 |z| = {e^{{\hat \rho}(z)} - 1\over e^{{\hat \rho}(z)} + 1},
\end{eqnarray}
where the right-hand side is monotone increasing in ${\hat \rho}(z)$.
We may apply a comparison theorem of Greene and Wu ([GW 1], Prop. 2.1, 
p. 26) to the two metrics, $ds^2$ and the hyperbolic metric $d{\hat s}^2$, 
on the disk $|z| \leq r$. The comparison
theorem states that for any smooth monotone 
increasing function $f$, one has 
\[ \bigtriangleup (f \circ \rho) \leq {\hat \bigtriangleup} (f \circ {\hat \rho}) \]
where $\rho$ and ${\hat \rho}$ are the distances to the origin in the 
metrics $ds^2$ and $d{\hat s}^2$ respectively, $\bigtriangleup$ and 
${\hat \bigtriangleup}$ are the Laplacians with respect to the two metrics, 
and the two sides are evaluated at points of the same level sets of the two 
metrics, i.e. $\rho = c$ on the left and ${\hat \rho} = c$ on the right, 
provided in two dimensions that the Gauss curvatures $K$ and 
${\hat K}$ satisfy $0 \ge K \ge {\hat K}$, with a similar 
condition on Ricci curvature in higher dimension. In our case we have 
$0 \ge K \ge -1 = {\hat K}$, and 
so we may apply the theorem. We note that the function 
\[ \log |z| = \log {e^{{\hat \rho}(z)} -1 \over e^{{\hat \rho}(z)} +1} \]
is harmonic with respect to $z$ and is therefore also harmonic 
with respect to any conformal metric on 
$0 < |z| < 1$. In other words, if we set 
\[ f(t) = \log {e^t -1 \over e^t + 1} \]
we have 
\[ {\hat \bigtriangleup} (f \circ {\hat \rho}) \equiv 0 \] 
for $0 < |z| < 1$. Since $f$ is monotone increasing, we may apply the 
Greene-Wu 
comparison theorem to conclude that 
\[ \bigtriangleup (f \circ \rho) \leq 0 \]
for $0 < |z| < r$, i.e. $f\circ \rho$ is superharmonic. 
For $z$ near $0$, we have 
$\rho(z) \sim \lambda(0) |z|$, and we may apply the minimum principle to the 
function 
\[ \log {1\over |z|} f(\rho(z)) = \log {1\over |z|}{e^{\rho(z)} -1 \over e^{\rho(z)} + 1}, \]
which is superharmonic in $0 < |z| < r$ and bounded near the origin, to 
conclude that it takes on its minimum on the boundary $|z| = r$. But since 
$\rho(z) \ge R$ on $|z| = r$, we have for $|z| < r$ that 
\[ \log {1\over |z|}{e^{\rho(z)} -1 \over e^{\rho(z)} + 1} \ge 
\log {1\over r}{e^R -1 \over e^R + 1} = 0,\]
by (2). Hence 
\[ {e^{\rho(z)} -1 \over e^{\rho(z)} + 1} \ge |z| = {e^{{\hat \rho}(z)} -1 \over e^{{\hat \rho}(z)} + 1}, \] 
by (3), which implies $\rho(z) \ge {\hat \rho}(z)$, proving the lemma.\qed

\bigskip
As an application of lemma 2.1, we have the following lemma:

\begin{lemma} 
Let $ds_n^2$ be a sequence of conformal metrics on the unit disk $D$ 
whose curvatures satisfy $-1 \leq K_n \leq 0$. Suppose that $D$ is a geodesic 
disk of radius $R_n$ with respect to the metric $ds_n^2$, where $R_n \rightarrow \infty$, and that the metrics $ds_n^2$ converge, uniformly on compact sets, 
to a metric $ds^2$. Then all distances to the origin with respect to $ds^2$ are greater than or equal to the corresponding hyperbolic distances in $D$. In 
particular, $ds^2$ is complete. 
\end{lemma}
\noindent{\sl Proof.}
For any point $z$ in $D$, let $\rho_n(z)$ be the distance from $0$ to $z$ in the metric $ds_n^2$ and let $\rho(z)$ be the distance in the limit metric $ds^2$. Let $|z| = r_n$ be the circle in $D$ of hyperbolic radius $R_n$. Explicitly, 
by Remark 2.1 above, 
\[ R_n = \log {1+r_n \over 1 - r_n}. \]
If we make the change of parameter $w = r_n z$, we may apply lemma 2.1 to 
the induced metric in $|w| < r_n$ and conclude that 
\[ \rho_n(z) \ge \log {1+|w| \over 1 - |w|} = \log {1+r_n |z| \over 1 - r_n |z|}. \]
As $n \rightarrow \infty$ we have $R_n \rightarrow \infty$ and $r_n \rightarrow 1$. Hence, by uniform convergence on compact sets, we have 
\[ \rho(z) = \lim_{n \rightarrow \infty} \rho_n(z) \ge \lim_{r_n \rightarrow 1} \log {1 + r_n |z| \over 1 - r_n |z|} = \log {1 + |z| \over 1 - |z|}, \]
which proves the lemma. \qed

\bigskip
{\bf Note}: Although we shall not make use of it, we remark that lemma 2.1 
also implies another dual form of the Ahlfors-Schwarz-Pick lemma, closer in form to the original:

\begin{lemma} Let $S$ be a simply-connected surface with a complete metric 
$ds^2$ whose Gauss curvature satisfies $-1 \leq K \leq 0$. 
If $S$ is mapped 
conformally onto the unit disc, then the distance between any two points 
of $S$ is greater than or equal to the hyperbolic distance between the 
corresponding points in the disk. 
\end{lemma}
\noindent{\sl Proof.}  Given two points $p, q$ of $S$, we may map $p$ onto the 
origin, and let $z$ be the image of the point $q$. Then the distance 
between $p$ and $q$ on $S$ is given by $\rho(z)$ in terms of the pull-back 
of the metric on $S$ onto the disk. For any $r$ such that $|z| < r < 1$, let 
${\hat \rho}(z)$ be the hyperbolic distance from $0$ to $z$ and let 
$\rho_r(w)$ be the pullback of the metric on $S$ to $|w| < r$ under the map 
$z = w/r$. Then, since $S$ is complete, we may apply lemma 2.1 to conclude 
that 
\[ {\hat \rho}(z) \leq \rho_r(w) = \rho_r(rz). \]
But as $r \rightarrow 1, \rho_r(rz) \rightarrow \rho(z)$, which proves the lemma. 
\qed

\bigskip
Note that Lemma 2.3 combined with the standard Ahlfors-Schwarz-Pick lemma 
implies a generalization of Ahlfors' lemma due to Yau ([Y1]; see also Troyanov
[T1]):  {\sl Let $S_1$ be a simply-connected Riemann surface with a complete metric 
$ds^2$ whose Gauss curvature satisfies $-1 \leq K \leq 0$, and let $S_2$ be a
Riemann 
surface with Gauss curvature bounded above by $-1$. Let $f: S_1 \rightarrow 
S_2$ be a holomorphic map. Then $f$ is distance decreasing}.

\bigskip
We also need the following more precise version of theorem 1.3; the proof 
follows exactly as in Ru[1]. 
\begin{theorem}(cf. Ru[1])
Let $x: M \rightarrow {\R }^m$ be a complete minimal surface immersed in 
${\R }^m$. Suppose that its generalized Gauss map $g$ omits the hyperplanes 
$H_1, \dots, H_q$  in ${\P }^{m-1}({\C })$ and $g(M)$ is contained in 
some ${\P }(V)$, where $V$ is a subspace of ${\C }^m$ of dimension $k$. 
Assume that $H_1 \cap {\P }(V), \dots, H_q \cap {\P }(V)$ are 
in general position in ${\P }(V)$ and $ q > k(k+1)/2$. Then
$g$ must be constant. 
\end{theorem}

\bigskip
The following theorem due to M. Green (see [G1]) shows that the 
complement of $2m+1$ hyperplanes in 
general position in ${\P }^m({\C })$ is complete Kobayashi hyperbolic. 
\begin{theorem} Let $H_1, \dots, H_q$ be hyperplanes in 
${\P }^m({\C })$, located in general position. If $q \ge 2m+1$, then 
$X = {\P }^m({\C }) - \cup_{j=1}^q H_j$ is complete hyperbolic and 
hyperbolically imbedded in ${\P }^m({\C })$. Hence, 
if $D \subset {\C }$ is the unit disc, and $\Phi$ is a subset of 
$Hol(D, X)$, then $\Phi$ is 
relatively locally compact in $Hol(D, {\P }^m({\C }))$, i.e. given a 
sequence $\{f_n\}$ in $\Phi$ there exists a subsequence which 
converges 
uniformly  on every compact subset of $D$ to an element of 
$Hol(D, {\P }^m({\C }))$. 
\end{theorem} 
For the notions of ``complete Kobayashi hyperbolicity'' 
and ``hyperbolically imbedded in 
${\P }^m({\C })$'', see Lang [L1]. 

\bigskip
Before going to the next section, 
we  recall here a standard definition. 
\begin{definition}
Let $f: M \rightarrow {\P }^n({\C })$ be a holomorphic map. 
Let $p \in M$. {\bf A local reduced representation of $f$ around $p$} is 
a holomorphic map ${\tilde f}: U \rightarrow {\C }^{n+1} - \{{\bf 0}\}$, 
such that 
${\P }({\tilde f}) = f$, where $U$ is a neighborhood of $p$, and ${\P }$ 
is the projection map of ${\C }^{n+1} - \{0\}$ onto ${\P }^n({\C })$.
\end{definition}

\section{Proof of the main theorem}
Let $x: M \rightarrow {\R }^m$ be a minimal surface, where 
$M$ is a connected, oriented, real-dimension 2 manifold without 
boundary, and $x=(x_1,\dots, x_m)$ is an immersion. Then $M$ is a Riemann surface
in the induced structure defined by local  isothermal 
coordinates $(u,v)$.
The generalized 
Gauss map of the minimal surface, $g=[{\partial x_1\over \partial 
z}: \dots: {\partial x_m \over \partial z}]: M \rightarrow Q_{m-2}({\C }) \subset {\P }^{m-1}({\C })$ is a holomorphic map, where $z = u+iv$. 
The metric $ds^2$ on $M$, induced from the standard metric in ${\R }^m$, 
is $ds^2 = \sum_{j=1}^m |{\partial x_j\over \partial z}|^2 dz d{\bar z}$, 
and the Gauss curvature $K$ is given by ([HO1] p.37)
\begin{eqnarray}
 K = - 4{|{\tilde g}\wedge {\tilde g}'|^2\over |{\tilde g}|^6} = -
4{\sum_{ j<k} |g_j g'_k -
g_k g'_j|^2 \over (\sum_{j=1}^m |g_j|^2)^3}
\end{eqnarray}
where ${\tilde g} = (g_1, \dots, g_m), g_j = {\partial x_j \over \partial z}, 1 \leq j \leq m.$ 

\bigskip
We will need the following lemma:
\begin{lemma}
Let $M$ be a Riemann surface. Let $f_n: M \rightarrow {\P }^m({\C })$ 
be a sequence of holomorphic maps converging uniformly on every compact subset of $M$ to a holomorphic map $f: 
M \rightarrow {\P }^m({\C })$. Given 
 ${\bf a}, {\bf b} \in {\P }^m({\C }^*)$,
  let $f_{{\bf a}, {\bf b}}$ 
be the meromorphic function (called coordinate 
function) defined by 
\[f_{{\bf a}, {\bf b}}|_U  = { \alpha({\tilde f}) \over  \beta({\tilde f})},\]
where 
${\tilde f}$
is a reduced representation of $f$ on $U$, and $\alpha, \beta \in {\C }^{{m+1}^*}$ 
such that  ${\bf a} = {\P }(\alpha), 
{\bf b} = {\P }(\beta)$. 
Assume that 
 $\beta({\tilde f}) \not\equiv 0$ on some $U$ (i.e. the image of $f$ 
is not contained in the hyperplane defined by ${\bf b}$). Let 
$p \in M$ be such that  $\beta({\tilde f})(p) \not= 0$, and $U_p$ be a neighborhood of $p$ such that $\beta({\tilde f})(z) \not= 0$ for $z \in U_p$;  then  
$\{f_{n_{{\bf a}, {\bf b}}}$\} converges uniformly on  $U_p$ to the meromorphic function 
$f_{{\bf a}, {\bf b}}$. 
\end{lemma}
Proof. Since  the image of $f$ 
is not contained in the hyperplane defined by ${\bf b}$,  the image of
$f_n$ is also 
not contained in the hyperplane defined by ${\bf b}$ for $n$ large enough. Since ${{\bf a}({\bf x})\over {\bf b}({\bf x})}$ is a rational function 
on ${\P }^m({\C })$ and $f_n$  converges uniformly on every compact subset of $M$
to $f$, the composition functions 
also converge compactly. 
This concludes the proof. \qed

\bigskip
\begin{lemma} Let $x^{(n)}= 
(x_1^{(n)}, \dots, x_m^{(n)}): M \rightarrow {\R }^m$ be a sequence 
of  minimal immersions, and $g^{(n)}: M \rightarrow Q_{m-2}({\C }) 
\subset {\P }^{m-1}({\C })$  the sequence of their 
(generalized) Gauss maps. Suppose that $\{g^{(n)}\}$ converges 
uniformly on every compact subset of 
$M$ to 
a non-constant 
holomorphic map $g: M \rightarrow  Q_{m-2}({\C }) 
\subset {\P }^{m-1}({\C })$
and that there is some $p_0 \in M$ such that for each $j, 1 \leq j \leq m$, $\{x_j^{(n)}(p_0)\}$ converges. Assume also that $\{|K_n|\}$ is uniformly 
bounded, 
where $K_n$ is the Gauss curvature of the minimal surface $x^{(n)}$. 
Then 

\bigskip 
(i) either a subsequence $\{K_{n'}\}$ of $\{K_n\}$ converges to zero or

(ii) a subsequence $\{x^{(n')}\}$ of $\{x^{(n)}\}$ converges to a 
minimal immersion, $x: M \rightarrow {\R }^m$, whose Gauss map 
is $g$.
\end{lemma}

Proof.  By assumption, $g$ is not constant and we may 
assume that
$|K_n| \leq 1$ in $M$, 
for each $n \in {\bf N}$. For every 
point $p \in M$ let $(U_p, z)$ be a complex local coordinate centered 
at $p$. Let 
${\tilde g^{(n)}}= (g^{(n)}_1, \dots, g^{(n)}_m)$ 
 where $g^{(n)}_i = 
{\partial x_i^{(n)} \over \partial z}, 1 \leq i \leq m$, and 
let ${\tilde g}= (g_1, \dots, g_m)$ be a local reduced 
representation of $g$ on $U_p$.
Since some $g_i(z)$ is non-zero for each $z$, we know that
$g(M)$ is not contained in some coordinate hyperplane. 
Without loss of generality, we 
assume that $g(M)$ is not contained in 
the first coordinate hyperplane $H_1 = \{[y_1:\dots: y_m] \in {\P }^{m-1}({\C })| y_1 = 0\}$.
 Let 
\[ M_1 = \{ p \in M| g(p) \not\in H_1, {\tilde g}(p)\wedge {\tilde g}'(p) 
\not= 0 \}. \]
Note that $M - M_1$ is a discrete set: namely, it consists of the zeros of 
$g_1$ (which are isolated, since $g(M) \not\subset H_1$, which is equivalent 
to $g_1 \not\equiv 0$) together with the common zeros of the components 
of ${\tilde g}\wedge {\tilde g}'$,
 which are the holomorphic functions
 $g_j g_k' - g_k g_j'$. In particular, 
\[ g_1 g_k' - g_k g_1'
 = g_1^2({g_k\over g_1})'\]
so that  ${\tilde g}\wedge {\tilde g}' \equiv 0$ implies that 
$g_k/g_1 = c_k$, a  constant for each $k$, so that ${\tilde g} = g_1(1, c_2, \dots, c_m)$ and the map $g$ would be constant, contrary to
assumption. Thus, the zeros of  ${\tilde g}\wedge {\tilde g}'$ are isolated and the points of $M - M_1$ are also isolated.

\bigskip
Let $p \in M_1$. 
Since $g(p) \not\in H_1$, there is a neighborhood $U_p$ of $p$ such that $g_1(z) \not\in H_1$, 
and $g^{(n)}(z) \not\in H_1$ for $n$ large enough, for every $z \in U_p$. 
Choosing $U_p$ sufficiently
small, we have that  $g_2/g_1, \dots, g_n/g_1$ are holomorphic and 
\[ 4{|{\tilde g}\wedge {\tilde g}'|^2/|g_1|^4\over (1+ \sum_{j=2}^m |g_j/g_1|^2)^3} = 4{\sum_{ j<k} |{g_j\over g_1} ({g_k\over g_1})' -
{g_k\over g_1} ({g_j\over g_1})'|^2 \over (1+ \sum_{j=2}^m |g_j/g_1|^2)^3}
\ge 2c_1,\]
in $U_p$, where $c_1$ is some positive constant.
Since $g^{(n)} \rightarrow g$ uniformly, by lemma 3.1, $\{g^{(n)}_j/g_1^{(n)}\}$ converges 
uniformly to $g_j/g_1$ on $U_p$, $1 \leq j \leq m$.
So we have 
\[  4{\sum_{ j<k} |{g^{(n)}_j\over g^{(n)}_1} ({g^{(n)}_k\over g^{(n)}_1})' -
{g^{(n)}_k\over g^{(n)}_1} ({g^{(n)}_j\over g^{(n)}_1})'|^2 \over (1+ \sum_{j=2}^m |g^{(n)}_j/g^{(n)}_1|^2)^3}
\ge c_1,\]
in $U_p$
and so, by (4),
\[ {c_1 \over |g^{(n)}_1|^2} 
\leq  4{\sum_{ l<k} |{g^{(n)}_l\over g^{(n)}_1} ({g^{(n)}_k\over g^{(n)}_1})' -
{g^{(n)}_k\over g^{(n)}_1} ({g^{(n)}_l\over g^{(n)}_1})'|^2 \over |g^{(n)}_1|^2(1+ \sum_{j=2}^m |g^{(n)}_j/g_1^{(n)}|^2)^3} = |K_n| \leq 1,\]
in $U_p.$
Therefore 
\[ c_1 \leq |g^{(n)}_1|^2\]
in $U_p$, for large $n$. 
Then
$\{g^{(n)}_1\}$ is relatively 
compact in ${\mathcal M}(U_p)$. Noticing that $M - M_1$ is discrete, 
by taking a subsequence, if necessary,  we can 
assume that the globally defined holomorphic 1-forms $\{g^{(n)}_1 dz\}$ 
converge  
on $M_1$, to a  holomorphic 1-form $h_1dz$ or to infinity, uniformly on every
compact subset of $M_1$. 
We consider  each case below:

\bigskip
Case 1. $\{g^{(n)}_1 dz\}$ converges to infinity  uniformly on every 
compact subset of $M_1$.

\bigskip
For $p \in M_1$, we have, by (4), 
\begin{eqnarray} 
 K_n(p) = -4{\sum_{ j<k} |{g^{(n)}_j\over g^{(n)}_1} ({g^{(n)}_k\over g^{(n)}_1})' -
{g^{(n)}_k\over 
g^{(n)}_1} ({g^{(n)}_j\over g^{(n)}_1})'|^2 \over |g^{(n)}_1|^2(1+ \sum_{j=2}^m |g^{(n)}_j/g^{(n)}_1|^2)^3} \rightarrow 0.
\end{eqnarray}

Let $p$ be a point such that 
$p\not\in M_1$ but also
$g(p) \not\in H_1$; then in a small disc of $U_p$, 
$D(2\epsilon)$, $g^{(n)}(z) \not\in H_1$ for $n$ large enough, $z \in D(2\epsilon)$. This means that $g^{(n)}_1$ is non-vanishing on $D(2\epsilon)$ and 
 $g^{(n)}_1$  converges to infinity on $\partial D(\epsilon)$. 
From the maximum principle 
we conclude that $\{g^{(n)}_1\}$ converges to  
infinity on $D(\epsilon)$.
Therefore we again have $K_n(p) \rightarrow 0$ by (4). 

Finally suppose that $g(p) \in H_1$, i.e. $g_1(p) = 0$. Since $g(p)$ is not 
contained in some coordinate hyperplane, 
we assume that $g(p) \not\in H_2$, where 
$H_2$ is the second coordinate hyperplane, 
$H_2 = \{[y_1:\dots: y_n] \in {\P }^{n-1}({\C })| y_2 = 0\}$. Therefore, 
on a small disc, $D(2\epsilon)$, $g^{(n)}(z) \not\in H_2$ for $n$ large enough, i.e. $g^{(n)}_2(z) \not= 0$, for $z \in D(2\epsilon)$, and $g^{(n)}_1$, 
$g_1$ have no zeros 
other than the point $p$ on $D(2\epsilon)$ for $n$ large enough.
By lemma 3.1, $\{{g^{(n)}_2\over g^{(n)}_1}\}$, as a  
sequence  of  non-vanishing
holomorphic functions,   converges uniformly on $\partial D(\epsilon)$. 
Clearly,  $\{{g^{(n)}_2\over g^{(n)}_1} 
g^{(n)}_1\}$  converges uniformly to infinity on $\partial D(\epsilon)$,
and therefore $g^{(n)}_2$ converges uniformly to infinity on $\partial
 D(\epsilon)$.
 Again from 
the maximum  principle, we conclude that $g^{(n)}_2$ converges 
to infinity on $D(\epsilon)$.  By (4), noticing that 
\[ |{\tilde g}^{(n)}\wedge {\tilde g}^{(n)\prime}|^2 /|g_2^{(n)}|^4 = \sum_{j < k} |{g^{(n)}_j\over g^{(n)}_2} ({g^{(n)}_k\over g^{(n)}_2})' -
{g^{(n)}_k\over g^{(n)}_2} ({g^{(n)}_j\over g^{(n)}_2})'|^2,\]
we have 
\[ K_n(p) = -4{\sum_{ j<k} |{g^{(n)}_j\over g^{(n)}_2} ({g^{(n)}_k\over
 g^{(n)}_2})' -
{g^{(n)}_k\over g^{(n)}_2} ({g^{(n)}_j\over g^{(n)}_2})'|^2 \over |g^{(n)}_2|^2(\sum_{j=1}^n |g^{(n)}_j/g^{(n)}_2|^2)^3} \rightarrow 0.\]

\bigskip
Thus, we have proved that $K_n(p) \rightarrow 0$ for all $p \in M$. 
This corresponds to  
case (i) of the lemma.

\bigskip
Case 2. $\{g^{(n)}_1 dz\}$ converges to a holomorphic 1-form, 
$h_1dz$, on $M_1$.

\bigskip 
Let $p \in M - M_1$. If $D(2\epsilon)$ is a small disc contained in $U_p$, 
as $\{g^{(n)}_1 \} \rightarrow h_1$ uniformly on 
$\partial D(\epsilon)$ and $g_1^{(n)}$ are holomorphic,
 using the maximum principle, we see that $\{g^{(n)}_1 \}$ is 
relatively compact on $D(\epsilon)$. Therefore
$h_1dz$ extends to a holomorphic 1-form on  $M$ and  
the global 1-forms $\{g^{(n)}_1 dz\}$ converge to $h_1dz$ on $M$.

\bigskip
We now prove that, for every integer $j, 2 \leq j \leq m$,  
the global 1-forms $\{g^{(n)}_j dz\}$ converge to a holomorphic form 
$h_jdz$ on $M$.  
Let $p \in M$ such that  $g(p) \not\in H_1$; then 
there is a neighborhood 
$U_p$ of $p$ such that $g_1$, $g_1^{(n)}$ have no zeros for $n$ large enough, 
$z \in U_p$. 
Since $g^{(n)}_j = {g^{(n)}_j\over g^{(n)}_1} g^{(n)}_1$, 
and by lemma 3.1, 
$\{{g^{(n)}_j\over g^{(n)}_1}\}$ 
converges uniformly on $U_p$, 
and $g^{(n)}_1$ also converges uniformly on $U_p$,  $\{g^{(n)}_j\}$ must
 converge uniformly on  $U_p$.
For the points $p$ such that $g(p) \in H_1$, if $D(2\epsilon) \subset U_p$ is 
small enough so that $g_1, g_1^{(n)}$ have no zeros other than $p$ for $n$ large enough on $D(2\epsilon)$, then
we just proved that $\{g^{(n)}_j\}$ is uniformly 
 convergent  on $\partial D(\epsilon)$. Since $g^{(n)}_j$ are holomorphic, 
 by the maximum principle, 
we have that $g^{(n)}_j$  converges uniformly on $D(\epsilon)$. 
Therefore the globally defined holomorphic 1-forms $g^{(n)}_jdz$ converge to 
$\omega_j = h_j dz$, $1 \leq j \leq m$. Obviously these $\omega_j$ satisfy
 the conditions in Theorem 2.1
 since the $g^{(n)}_j$ satisfy the conditions in Theorem 2.1 for every $n$.
 So they define a minimal surface $x: M \rightarrow 
{\R }^m$ whose Gauss map is $g$. 
\qed

\bigskip
We now prove the main theorem.

\bigskip
\noindent{\sl Proof of the Main Theorem}.
 
 Suppose the theorem is not true. We will construct a nonflat
 complete minimal 
surface whose Gauss map omits a set of
 hyperplanes in general position, thus getting 
 a contradiction with theorem 2.3. 
So suppose  the conclusion of the theorem is not true; then there is 
a sequence of (non complete) minimal surfaces $x^{(n)}: 
M_n \rightarrow {\R }^m$ and points $p_n \in M_n$ 
such that $|K_n(p_n)|d_n^2(p_n) \rightarrow \infty$, and such that the Gauss map $g^{(n)}$ of $x^{(n)}$ omits a fixed set of $q$ hyperplanes in general position, with $q > {m(m+1)/2}$.

\bigskip 
We claim that the surfaces $M_n$ can be chosen so that 
\begin{eqnarray}
K_n(p_n) = -1,~~ -4 \leq K_n \leq 0 ~on ~M_n ~for~ all~ n 
~and~ d_n(p_n) \rightarrow \infty.
\end{eqnarray}
 We now prove the claim. Without 
loss of generality, we can assume that $M_n$ is a geodesic disk centered at $p_n$. Let $M_n' = \{p\in M_n: d_n(p, p_n) 
\leq d_n(p_n)/2\}$. Then $K_n$ is uniformly bounded on $M_n'$ 
and $d_n'(p)=$ distance of $p$ to the boundary of $M_n'$ 
tends to zero as $p \rightarrow \partial M_n'$. Hence
$|K_n(p)|(d_n'(p))^2$ has a maximum at a point $p_n'$ 
interior to $M_n'$. Therefore
\[ |K_n(p_n')|d'_n(p_n')^2 \ge |K_n(p_n)|d'_n(p_n)^2 =
{1\over 4}|K_n(p_n)|d_n^2(p_n) \rightarrow \infty. \]
So we can replace the $M_n$ by the $M_n'$, with $|K_n(p_n')|
d'_n(p_n')^2 \rightarrow \infty.$  We rescale $M_n'$ to make $K_n(p_n') = -1$.
By the invariance under  
scaling of the quantity $K(p)d(p)^2,$ we will have $d_n'(p_n') \rightarrow 
\infty$; here, without causing confusion, 
we use the same notation $d_n'$ to denote 
the geodesic distance with respect to the rescaled metric. 
Again we can assume that $M_n'$ is a geodesic disc centered at $p_n'$, and let 
\[ M_n'' = \{p \in M_n'| d_n(p, p_n') < {d_n'(p_n')\over 2}\}.\]
Then $p \in M_n''$ implies that $d_n'(p) \ge {d_n'(p_n')\over 2}$
and 
\[ |K_n(p)|{d'_n(p_n')^2\over 4} \leq |K_n(p)|d'_n(p)^2 \leq 
|K_n(p_n')|d'_n(p_n')^2 = d'_n(p_n')^2. \]
Therefore
$|K_n(p)| \leq 4$ on $M_n''$. Furthermore, 
$d_n''(p_n') = d(p_n', \partial M_n'') = d_n'(p_n')/2 \rightarrow 
\infty$. This proves the claim. 

\bigskip
By translations of
 ${\R }^m$ we can assume that $x^{(n)}(p_n) = {\bf 0}$.
We can 
also assume that $M_n$ is simply connected, by taking its universal 
covering, if necessary. 
By the uniformization theorem, $M_n$ is conformally equivalent to 
either the unit disc $D$ or the complex plane ${\C }$, and 
we can suppose that $p_n$ maps onto $0$ for each $n$. 
But the  
case that $M_n$ is conformally equivalent to ${\C }$ is impossible because 
the condition that 
$g^{(n)}$ misses more than $m(m+1)/2$ hyperplanes in general 
position in ${\P }^{m-1}({\C })$, implies, by Picard's theorem, that 
$g^{(n)}$ is constant, so $K_n \equiv 0$, which 
contradicts 
the condition that $|K_n(0)| = 1$.  
So we have constructed a sequence of  minimal surfaces, 
$x^{(n)}: D \rightarrow {\R }^3$, satisfying (6).
 Since, by theorem 2.4,  
${\P }^{m-1}({\C })$ minus $2m-1$ hyperplanes
is complete 
Kobayashi hyperbolic, and $m(m+1)/2 \ge 2m-1$,  a subsequence of 
 generalized Gauss maps 
$g^{(n)}$ of $x^{(n)}$ exists$-$without of loss generality we 
assume $g^{(n)}$ itself$-$such that $g^{(n)}: D \rightarrow {\P }^{m-1}({\C })$
 converges uniformly on every compact subset of $D$
 to a map $g: D \rightarrow {\P }^{m-1}({\C })$. 
 
\bigskip
We now claim that $g$ is non-constant. Suppose not, i.e. $g$ is a constant 
map, and $g$ maps the disk $D$ onto a single point $P$.  
Let $H$ be any hyperplane not containing the point $P$, and let  $U$, $V$ be 
disjoint neighborhoods of $H$ and $P$ respectively. Let $C$ be the constant 
in theorem 1.2 such that 
\[ |K(p)|^{1/2}d(p) \leq C \]
for any minimal surface in ${\R }^m$ whose Gauss map omits the neighborhood 
$U$ of $H$, where $p$ is a point of $S$ and $d(p)$ is the geodesic distance of 
$p$ to the boundary of $S$. Choose $r < 1$ such that the hyperbolic distance 
$R$ of $z=0$ to $|z| = r$ satisfies $R > C$. Since $g^{(n)}$ converges 
uniformly to $g$ on $|z| \leq r$, the image of $|z| = r$ lies in the 
neighborhood $V$ of $q$ for sufficiently large $n$, say $n \ge n_0$. It 
follows that for $n \ge n_0$, the image of the disk $|z| \leq r$ 
under $g^{(n)}$ omits the neighborhood $U$ of $H$ and we may therefore apply the above 
inequality to conclude 
\[ |K_n(0)|^{1/2} d_n(r) \leq C\]
where $d_n(r)$ is the geodesic distance from the origin to the boundary of 
the surface $x^{(n)}: D(r) \rightarrow {\R }^m$. But $|K_n(0)| = 1$ for all 
$n$, and hence $d_n(r) \leq C$ for $n \ge n_0$. On the other hand, we get a 
lower bound for $d_n(r)$ from lemma 2.1. The surface $x^{(n)}: \{|z| < 1\}
\rightarrow {\R }^m$ is a geodesic disk of radius $R_n$. If we 
reparametrize by $w = r_n z$ where $|w| = r_n$ has hyperbolic radius $R_n$, 
then the circle $|z| = r$ corresponds to $|w| = r_n r$, and by lemma 2.1, 
the distance in the surface metric from the origin to any point on the 
circle $|z| = r$, or equivalently, $|w| = r_nr$, is greater than or 
equal to the hyperbolic distance from $0$ to $|w| = r_n r$. 
But as $n \rightarrow \infty$, $R_n 
\rightarrow \infty$ and $r_n \rightarrow 1$, so that the hyperbolic 
radius of $|w| = r_n r$ tends to the hyperbolic radius of $|w| = r$, 
which is $R$. Since by assumption $R > C$ we have for $n$ sufficiently large 
that the 
surface distance from $z=0$ to $|z| = r$ is greater than $C$, contradicting 
the earlier bound $d_n(r) \leq C$. Thus we conclude that 
the limit function $g$ 
can not be constant. 

\bigskip
Therefore the hypotheses of lemma 3.2 are satisfied. 
Since $|K_n(0)| = 1$, the possibility (i) of lemma 3.2 cannot happen. 
Thus, a subsequence $\{x^{(n')}\}$ of $\{x^{(n)}\}$ 
converges to a minimal immersion $x: D \rightarrow {\R }^m$, whose Gauss map is $g$.
By (6) and by lemma 2.2, $x$ is complete. 
By assumption,
$g^{(n)}$ omits hyperplanes $H_1, \dots, H_q$ in 
${\P }^{m-1}({\C })$, located in general position, $q > m(m+1)/2$. 
By Hurwitz's theorem(theorem 2.2), either $g$ omits these hyperplanes, or 
the image of $g$ lies in some of these hyperplanes. Say 
$g(M) \subset \cap_{j=1}^k H_{j} 
= {\P }(V)$, where $V$ is a subspace of ${\C }^m$ of dimension $m-k$, and
 $g(M): M \rightarrow {\P }(V)$ omits the hyperplanes $H_{k+1}\cap 
(\cap_{j=1}^k H_{j}), 
\dots, H_{q}\cap (\cap_{j=1}^k H_{j})$ in ${\P }(V)$.
Since the hyperplanes $H_{k+1}\cap 
(\cap_{j=1}^k H_{j}), 
\dots, H_{q}\cap (\cap_{j=1}^k H_{j})$ in ${\P }(V)$ are still in general position in ${\P }(V)$ because $H_1, \dots, H_q$ are in general position in ${\P }^{m-1}({\C })$, and $q-k >  m(m+1)/2 - k  \ge
(m-k)(m-k+1)/2$, 
it follows from  theorem 2.3 that $g$ is constant. But we have 
just proved that $g$ is not constant. This leads to a contradiction.
Therefore the main theorem is proved.  
\qed


\begin{thebibliography}{}

\bibitem
[A 1]{} L. V. Ahlfors, {\sl \"Uber die Anwendung differentialgeometrischer Methoden zur Untersuchung von \"Uberlagerungsfl\"achen,} Acta Soc. 
Sci. Fennicae Nova Ser. A 2, (1937), 1-17.   

\bibitem
[B 1]{} R. Brody, {\sl Compact manifolds and hyperbolicity}, 
Trans. Amer. Math. Soc., 235(1978), 213-219.

\bibitem
[CO 1]{} S. S. Chern and R. Osserman, {\sl Complete minimal surfaces in euclidean n-space}, J. d' Analyse Math., 19(1967), 15-34.

\bibitem
[F 1]{} H. Fujimoto, {\sl On the number of exceptional values of the Gauss map 
of minimal surfaces}, J. Math. Soc. Japan, 40(1988), 235-247.

\bibitem
[F 2]{} H. Fujimoto, {\sl Modified defect relations for the 
Gauss map of minimal surfaces, II}, J. Differential Geometry 31(1990), 365-385.

\bibitem
[F 3]{} H. Fujimoto, {\sl Value distribution theory of the Gauss map 
of minimal surfaces in ${\R }^m$}, Aspects of Math., E21(1992).


\bibitem
[G1]{} M. Green, {\sl The hyperbolicity of the complement of $2n+1$ hyperplanes 
in general position in ${\P }^n$ and related results}, Proc. Amer. Math. 
Soc., 66(1977), no. 1, 109-113. 

\bibitem
[GW1]{} R. E. Greene and H. Wu,  {Function theory on manifolds which possess 
a pole}, Lecture Notes in Mathematics, No. 699, Springer-Verlag, 1979.

\bibitem
[HO 1]{} D. Hoffman and R. Osserman, {\sl The geometry of the generalized Gauss
map}, Memoirs of the A.M.S., 236(1980). 

\bibitem
[L 1]{} S. Lang, {Introduction to complex hyperbolic space}, Springer-Verlag, 1987. 

\bibitem
[O 1]{} R. Osserman, {\sl Global properties of minimal surfaces in $E^3$ and $E^m$}, Ann. of Math., 80(1964), 340-364.

\bibitem
[O 2]{} R. Osserman, {\sl A survey of minimal surfaces}, 2nd edition, Dover Publ. Inc., New York, 1986. 

\bibitem
[Ro 1]{} A. Ros, {\sl The Gauss map of minimal surfaces}, Preprint. 


\bibitem
[R 1]{} M. Ru, {\sl On the Gauss map of minimal surfaces immersed in 
${\R }^n$}, J. Differential Geometry, 34(1991), 411-423.

\bibitem
[T 1]{} M. Troyanov, {\sl The Schwarz lemma for nonpositively curved 
Riemannian surfaces}, Manuscripta Math., 72(1991), 251-256.

\bibitem
[Y 1]{} S. T. Yau, {A general Schwarz lemma for K\"ahler manifolds}, 
Amer. Jour. Math., 100(1978), 197-203. 


\end{thebibliography}
\end{document}